# Prime Reciprocal Digit Frequencies and the Euler Zeta Function

Subhash Kak

1. The digit frequencies for primes are not all equal. The least significant digit for primes greater than 5 can only be 1, 3, 7, or 9. For binary representation of primes, the frequency of 1 is higher than that of 0.

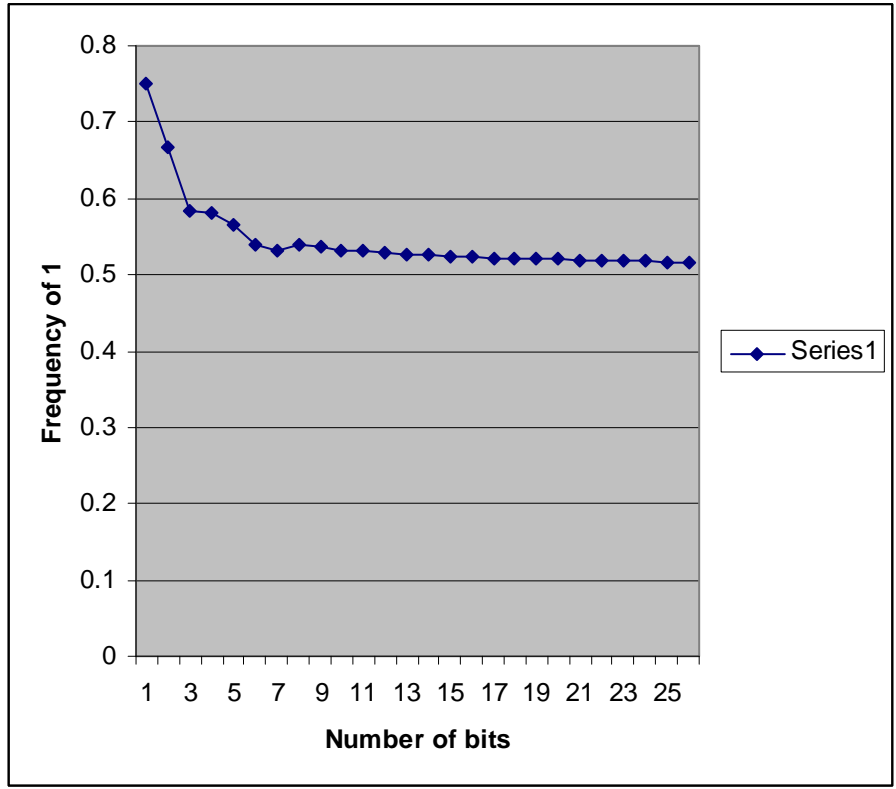

**Figure 1:** Frequency of 1 in a binary, constant bit representation of primes

The number of 0s and 1s for all primes of with respect to different binary lengths from 2 to 27 is given in Table 1. Thus for all primes of binary length 3, we have the primes 2, 3, 5, and 7 which in the binary form are 010, 011, 101, and 111, with four 0s and 8 1s. Likewise, for all primes of bit length 4, we count the primes 2, 3, 5, 7, 11, and 13 corresponding to the sequences 0010, 0011, 0101, 0111, 1011, 1101, which gives us ten 0s and 14 1s.

If one were to take a variable binary representation of primes, the frequency of 0s will naturally be less. Thus for 27-bit representation, the count of 0s in a variable bit representation is 90971028 and that of 1s is 106034897 (compare the constant bit values in Table 1).

**Table 1:** Frequency of 0 and 1 in a constant bit representation of primes

| Bits | Count 0 | Count 1 | Total 0+1 |
|---|---|---|---|
| 1 | | | |
| 2 | 1 | 3 | 4 |
| 3 | 4 | 8 | 12 |
| 4 | 10 | 14 | 24 |
| 5 | 23 | 32 | 55 |
| 6 | 47 | 61 | 108 |
| 7 | 100 | 117 | 217 |
| 8 | 202 | 230 | 432 |
| 9 | 403 | 470 | 873 |
| 10 | 798 | 922 | 1720 |
| 11 | 1592 | 1807 | 3399 |
| 12 | 3171 | 3597 | 6768 |
| 13 | 6293 | 7071 | 13364 |
| 14 | 12578 | 14022 | 26600 |
| 15 | 24987 | 27693 | 52680 |
| 16 | 49796 | 54876 | 104672 |
| 17 | 99190 | 109077 | 208267 |
| 18 | 197699 | 216301 | 414000 |
| 19 | 394227 | 430183 | 824410 |
| 20 | 785804 | 854696 | 1640500 |
| 21 | 1567419 | 1700412 | 3267831 |
| 22 | 3127966 | 3382868 | 6510834 |
| 23 | 6242519 | 6733230 | 12975749 |
| 24 | 12464093 | 13404811 | 25868904 |
| 25 | 24887586 | 26704639 | 51592225 |
| 26 | 49698098 | 53204936 | 102903034 |
| 27 | 99261034 | 106034897 | 205295931 |

The principal question of interest to us is whether the excess of 1s is counterbalanced by the higher frequency of 0 compared to 1 in the representation of prime reciprocals?

2. There is indeed evidence [1],[2] that for prime reciprocals in binary expansions the frequency of 0 is slightly larger than that of 1. For example, the number of cases where 0s exceed 1s compared to where 1s exceed 0s is in the proportion 7:1 in the primes in the range 50,000 to 60,000. Likewise, in the range 800,011 to 999,983, the number of cases where 0s exceed 1s in the prime reciprocals is 3609, whereas the number of cases where 1s exceed 0s is only 641 (the number of cases where 0s and 1s were equal was 10297, which are the maximum-length cases).

3. It is possible that the frequency of 0 being slightly larger than that of 1 does not hold up as the number of primes tested increases much beyond 999,983, which was the limit of the most recent experiment. But it is more likely that the difference will continue to hold up, especially if we see this excess as counterbalancing the excess of 1s in the binary representation of primes.

4. One would like to know how the difference between 0s and 1s changes as the range increases. With this in mind, it may be well worth studying the possible relationship of the reciprocal frequencies with the Euler and the Riemann zeta functions.

5. The Euler zeta function ζ(s) is given by:

$$\zeta(s) = \sum_{n \geq 1} 1/n^s$$

$$= \frac{1}{1^s} + \frac{1}{2^s} + \frac{1}{3^s} + \frac{1}{4^s} + \ldots$$

$$\zeta(1) = 1 + \frac{1}{2} + \frac{1}{3} + \frac{1}{4} + \ldots$$

This is the sum of all reciprocals, which is the harmonic series whose sum is infinity (the proof of which is elementary).

6. By the Euler identity

$$\zeta(s) = \prod_{p=prime} \frac{1}{1 - p^{-s}}$$

Therefore,

$$\zeta(1) = \prod_{p=prime} \frac{1}{1 - p^{-s}}$$

$$= \prod_{p=prime} \frac{p}{p-1}$$

Or,

$$\frac{1}{\zeta(1)} = \prod_{p=prime} \frac{p-1}{p}$$

Since $\frac{p-1}{p}$ is merely a cyclic displacement [3],[4] of the digits of the reciprocal expansion of 1/p, it is clear that there is some relationship between the Euler zeta function and the prime reciprocals. But it is not a direct relationship since it involves the multiplication of prime reciprocals. Perhaps this multiplication implies a frequency property related to all prime reciprocals.

$\frac{1}{\zeta(1)}$ is also written as a Dirichlet series $\sum_{n=1}^{\infty} \frac{\mu(n)}{n^s}$, where μ(n) is the Möbius function which is +1 if n is square-free with even number of distinct factors, -1 if it is square free with odd number of distinct factors, and 0 if it is not square free. Therefore,

$$\frac{1}{\zeta(1)} = 1 - 1/2 - 1/3 - 1/5 + 1/5 - 1/7 + 1/10 - 1/11 - 1/13 + 1/14 + 1/15 - \ldots$$

Considering binary prime reciprocal sequences, the difference between addition and subtraction may not matter in the computation of group frequencies.

7. Similarly,

$$\frac{1}{\zeta(2)} = \prod_{p=prime} \frac{p^2-1}{p^2}$$

This suggests that it would also be useful to study the expansions of reciprocals of prime powers. Since ζ(2) is $\pi^2/6$, could this imply a corresponding interesting property for such expansions?

8. By taking logs of both sides in the inverse expansion of the Euler zeta function,

$$\ln \frac{1}{\zeta(1)} = \sum_{p=prime} \ln \frac{(p-1)}{p}$$

Likewise, we have

$$\ln \frac{1}{\zeta(2)} = \sum_{p=prime} \ln \frac{(p^2-1)}{p^2}$$

9. It is possible that further insights into prime reciprocals may be obtained from s complex, or, in other words, from the Riemann zeta function.

10. It is convenient to generate binary expansions of the prime reciprocal 1/p by [4],[5]:

    $$a(i) = 2^i \bmod p \bmod 2$$

    It is worthwhile to determine if this phenomenon holds for non-binary cases. To generate decimal expansions of 1/p, one may use the following formula:

    A preliminary experiment indicates that in base 10 also the frequency of 0 is larger than that of the other digits.

11. Another interesting question to ask is what is the largest group of consecutive prime reciprocals that are non-maximum length in base 2? For example, the 6 consecutive prime reciprocals 970279, 970297, 970303, 970313, 970351, 970391 and 989647, 989663, 989671, 989687 989719, 989743 are all non-maximum length in base 2.

12. For some applications of prime reciprocal sequences, see [6]-[10].